\documentclass[12pt]{article}
\def\q{\hfill\rule{1ex}{1ex}}
\def\0{\emptyset}

\def\p{{\bf Proof.} \quad}
\def\q{\hfill\rule{1ex}{1ex}}

\def\n{\noindent}

\usepackage{graphicx}

\begin{document}
\title{\bf Matching preclusion and strong matching preclusion of the bubble-sort star graphs}
\author{{\small\bf
Xin Wang$^1$}\thanks{email: wangxin01@hnu.edu.cn}\quad\quad
{\small\bf Chaoqun Ma$^1$}\quad\quad
{\small\bf Jia Guo$^2$}\thanks{email: guojia199011@163.com}
\\
{\small $^1$Business School, Hunan University, Changsha, China, 410082}\\
{\small $^2$School of Software, Northwestern Polytechnical University, Xi'an, Shaanxi, China, 710072}
}

\date{}
\maketitle\baselineskip 15.5pt

\begin{abstract}
\baselineskip=0.5cm
Since a plurality of processors in a distributed computer system working in parallel, to ensure the fault tolerance and stability of the network is an
important issue in distributed systems. As the topology of the distributed network can be modeled as a graph, the (strong) matching preclusion in graph theory can be used as a robustness measure for missing edges in parallel and distributed networks, which is defined as the minimum number of (vertices and) edges whose deletion results in the remaining network that has neither a perfect matching nor an almost-perfect matching. The bubble-sort star graph is one of the validly discussed interconnection networks related to the distributed systems.
In this paper, we show that the strong matching preclusion number of an $n$-dimensional bubble-sort star graph $BS_n$ is $2$ for $n\geq3$ and each optimal strong matching preclusion set of $BS_n$ is a set of two vertices from the same bipartition set.
Moreover, we show that the matching preclusion number of $BS_n$ is $2n-3$ for $n\geq3$ and that every optimal matching preclusion set of $BS_n$ is trivial.

\vskip 0.3cm


{\bf Keywords:} Distributed systems, Network topology, Matching preclusion, Fault tolerance and stability

\end{abstract}

\vskip.3cm

\newpage
\n{\large\bf 1.\quad Introduction}

With the development of distributed computer technology, the importance of the interconnection network has become increasingly prominent. However, communication failures among multiple processors in the distributed network are inevitable, which has led to research on the fault tolerance and stability of the interconnection network becoming a current hot spot. The fault tolerance and stability are crucial factors in establishing and optimizing network topology. For example, the blockchain is a typical distributed network consisting of multiple processors, and the fault tolerance and stability of the network determine its security and availability. If a faulty vertex in a network is matched by a special matching, then the tasks running on the faulty vertex are able to be shifted to another vertex by the matching in the event of (vertex and) edge failure, which can enhance the fault tolerance and stability of the interconnection network. Thus the (strong) matching preclusion as a measure of the performance of the interconnection networks has been proposed and studied in recent years.

The interconnection network can be represented by an undirected graph $G$. The (strong) matching preclusion number of a graph $G$ is the minimum number of (vertices and) edges in $G$ whose deletion results in the remaining graph with neither perfect matchings nor almost-perfect matchings.
The concept of matching preclusion was firstly introduced as a measure of robustness of networks in the case of edge failure by Brigham et al. in [1]. To extend this concept, Park et al. in [14] gave the concept of strong matching preclusion.
The (strong) matching preclusion of many well-known interconnection networks has been studied, such as the star graphs [4], the augmented cubes [6],
the $k$-ary $n$-cubes [18]. Many other results can be seen in [3], [5], [10]-[13].

In this paper, we deal with the strong matching preclusion and the matching preclusion of the bubble-sort star graph $BS_n$, which has many nice properties. $BS_n$ is edge-bipancyclic, vertex-bipancyclic and bipancyclic for $n\geq3$ [9]. Zhao et al. [23] gave the generalized connectivity of $BS_n$ and Wang et al. [19] studied the diagnosability of $BS_n$ with missing edges. Many other properties of $BS_n$ have been investigated [8], [17], [20]-[22], [24]. To our knowledge, the research results of this paper are innovative and theoretically meaningful for distributed systems.

The remainder of this paper is organized as follows:
In section 2, we introduce some definitions and give the strong matching preclusion number of $BS_n$. We consider the matching preclusion of $BS_n$ in section 3. In section 4, we make a conclusion.

\vskip.3cm

\n{\large\bf 2.\quad Preliminaries}

Let $G$ be a simple connected graph with vertex set $V(G)$ and edge set $E(G)$. The {\em order} of a graph $G$, denoted by $|V(G)|$, is the number of vertices in $G$.
The {\em degree} of a vertex $v\in V(G)$, denoted by $d_G(v)$, is the number of edges that are incident with $v$ in $G$. A graph $G$ is {\em $k$-regular} if $d_G(v)=k$ for every $v\in V(G)$. Let $e_1=(v_1,v_2)$, $e_2=(u_1,u_2)$ be any two edges in a graph $G$. If $\{v_1,v_2\}\cap\{u_1,u_2\}=\emptyset$, then edges $e_1$ and $e_2$ are {\em independent} in $G$.
A graph $G$ is {\em bipartite} if there exist two non-empty bipartition sets $V_1,V_2\subseteq V(G)$ with $V_2=V(G)-V_1$ such that every edge of $G$ has one end vertex in $V_1$ and the other in $V_2$. The graph $K_{1,n-1}$ is a bipartite graph with vertex set $V(K_{1,n-1})=V_1\cup V_2=\{u_1,u_2,\cdots,u_n\}$ and edge set $E(K_{1,n-1})=\{(u_1,u_i)~|~i=2,3,\cdots,n\}$, where $V_1=\{u_1\}$, $V_2=\{u_2,u_3,\cdots,u_n\}$. The graph $P_n$ is a {\em path} with vertex set $V(P_n)=\{u_1,u_2,\cdots,u_n\}$ and edge set $E(P_n)=\{(u_i,u_{i+1})~|~i=1,2,\cdots,n-1\}$. Let $G,H$ be two graphs, then $G\cup H$ is a graph with vertex set $V(G)\cup V(H)$ and edge set $E(G)\cup E(H)$.

Let $G$ be a connected graph and $M\subseteq E(G)$. Then $M$ is a {\em matching} of $G$ if every vertex in $G$ is incident with at most one edge in $M$. The vertex which is incident with one edge in a matching $M$ is {\em matched} by $M$; vertices which are not incident with any edge of $M$ are called {\em unmatched} by $M$.
A matching $M\subseteq E(G)$ is a {\em perfect matching} if every vertex in $G$ is matched by $M$. A matching $M\subseteq E(G)$ is an {\em almost-perfect matching} if every vertex except one in $G$ is matched by $M$, while the exceptional vertex is unmatched by $M$. The order of $G$ is even if $G$ has a perfect matching, and the order of $G$ is odd if $G$ has an almost-perfect matching.

Let $F^*\subseteq V(G)\cup E(G)$, $F^*_E=\{(u,v)\in E(G)~|~ \{u,v\}\cap( V(G)\cap F^*)\neq \emptyset\}$ and $G-F^*$ be the subgraph of $G$ with vertex set $V(G)-F^*\cap V(G)$ and edge set $E(G)-F^*\cap E(G)-F^*_E$. The set $F^*$ is a {\em strong matching preclusion set} of $G$, if $G-F^*$ has neither perfect matchings nor almost-perfect matchings and the minimum size of such set $F^*$, denoted by $smp(G)$, is the {\em strong matching preclusion number} of $G$. A strong matching preclusion
set $F^*$ of graph $G$ is {\em optimal} if $|F^*|=smp(G)$.
Let $F\subseteq E(G)$ and $G-F$ be the subgraph of $G$ with vertex set $V(G)$ and edge set $E(G)-F$.
The {\em matching preclusion number} of $G$, denoted by $mp(G)$, is the minimum size of $F$ such that $G-F$ has no perfect matchings or almost-perfect matchings and such edge set $F$ is an {\em optimal matching preclusion set} of $G$. An optimal matching preclusion set $F$ of graph $G$ is {\em trivial} if every edge of $F$ is incident with exactly the same vertex in $G$. Obviously, $mp(G)=0$ if graph $G$ has neither perfect matchings nor almost-perfect matchings.

Let $n_1,n_2$ be two integers with $n_2> n_1\geq1$ and $[n_1,n_2]=\{n_1,n_1+1,\cdots,n_2-1,n_2\}$. Let $n,i,j$ be three integers with $n\geq2$ and $i,j\in[1,n]$. A permutation $x_1x_2\cdots x_n$ on $[1,n]$ is a {\em transposition} if $x_i=j,~x_j=i$ and $x_k=k$ for every $k\in[1,n]-\{i,j\}$. For any transposition $x_1x_2\cdots x_n$ with $x_i=j$ and $x_j=i$, set $\langle i,j\rangle=x_1x_2\cdots x_n$ for simplicity.
Let $a=a_1a_2\cdots a_n$, $b=b_1b_2\cdots b_n$ be two permutations on $[1,n]$ and ``$\circ$" be an operation on a transposition $\langle i,j\rangle$ such that $a=b\circ \langle i,j\rangle$ if and only if $a_i=b_j$, $a_j=b_i$ and $a_k=b_k$ for every $k\in [1,n]-\{i,j\}$.

\vskip.2cm
{\bf Definition 2.1 [7]} {\em Let $n$ be an integer with $n\geq2$ and $S$ be the set containing all permutations on $[1,n]$. An $n$-dimensional bubble-sort star graph $BS_n$ has vertex set $V(BS_n)=S$. For any two distinct vertices $u$ and $v$ in $V\left(BS_n\right)$, $u$ is adjacent to $v$ if and only if one of the following conditions holds:

$(1)$ $u=v\circ\langle1,i\rangle$ for some $i\in\left[2,n\right]$;

$(2)$ $u=v\circ\langle i-1,i\rangle$ for some $i\in\left[3,n\right]$.
}

\vskip.2cm

By Definition 2.1, $BS_n$ is a bipartite and $(2n-3)$-regular graph with $n!$ vertices.
For every $i\in[1,n]$, let $BS_n^i$ be the induced subgraph of $BS_n$ with $V(BS_n^i)=\{a=a_1a_2\cdots a_n\in V(BS_n)~|~a_n=i\}$.
Thus $BS_n^i\cong BS_{n-1}$ for every $i\in[1,n]$. For a subset $I\subseteq[1,n]$, let $BS_n^I$ be the subgraph of $BS_n$ induced by the vertex set $\cup_{i\in I}V(BS_n^i)$.
We illustrate $BS_2$, $BS_3$ and $BS_4$ in Fig. 1.

\begin{figure}[ht]
   \begin{center}
    \includegraphics[scale=1]{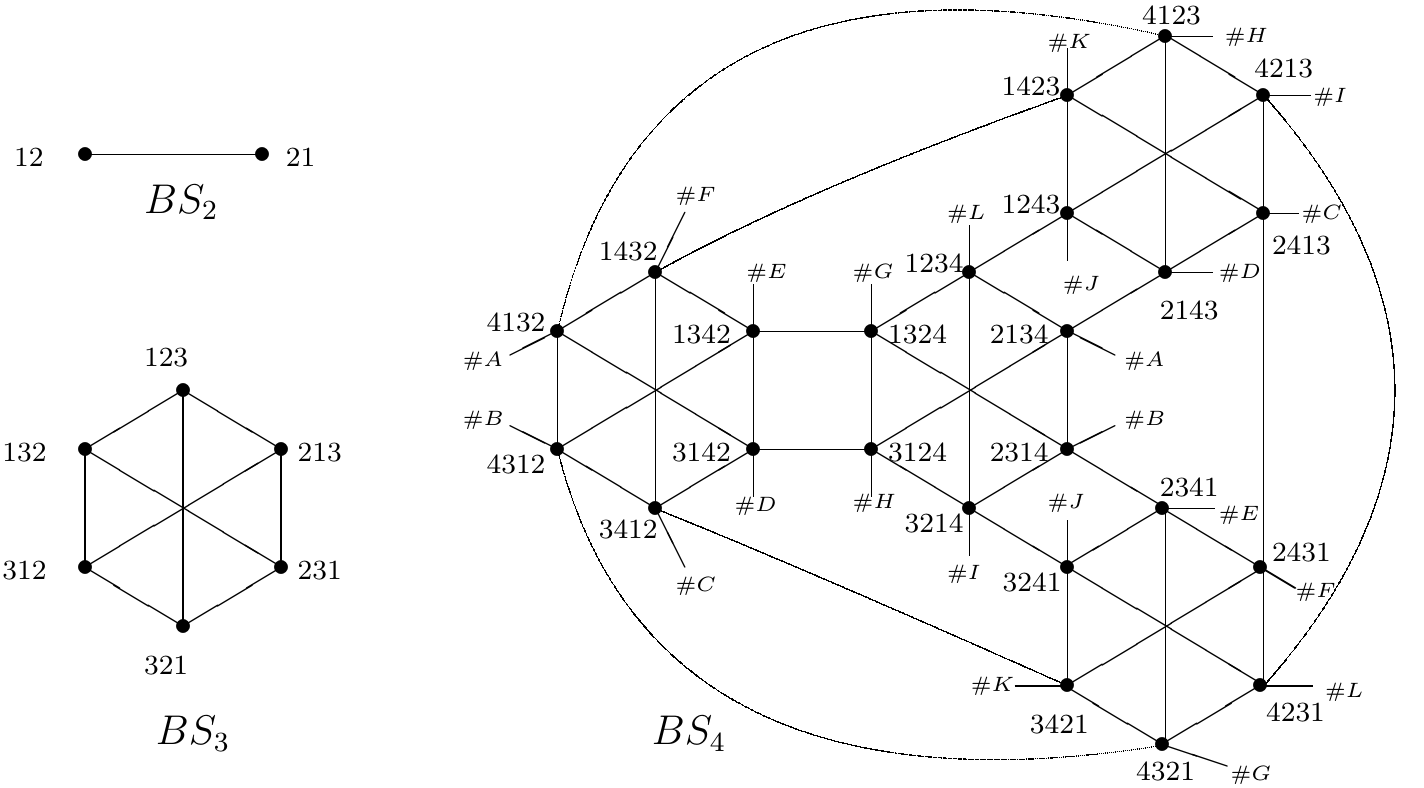}
    \end{center}
    \vspace{-0.5cm}\caption{\label{F2-5} $BS_n$ for $n=2,3,4.$}
\end{figure}

\vskip.2cm
Now we give some properties of $BS_n$.

\vskip.2cm
{\bf Lemma 2.2 [2]} {\em $BS_n$ is vertex transitive for $n\geq2$.}

\vskip.2cm

{\bf Lemma 2.3 [20]} {\em $BS_n$ is a special Cayley graph for $n\geq2$.}

\vskip.2cm

{\bf Lemma 2.4 [2]} {\em There exist $2(n-2)!$ independent edges connecting $BS_n^i$ and $BS_n^j$ for any $i,j\in[1,n]$ with $i\neq j$ and $n\geq3$.}

\vskip.2cm
{\bf Lemma 2.5 [2]} {\em For any two distinct vertices $u,v\in V(BS_n^i)$ with some $i\in[1,n]$ and $n\geq3$, $\{u\circ\langle1,n\rangle,~u\circ\langle n-1,n\rangle\} \cap \{v\circ\langle1,n\rangle,~v\circ\langle n-1,n\rangle\}=\emptyset$.}
\vskip.2cm

{\bf Lemma 2.6 [16]} {\em Let $G$ be a connected $k$-regular and bipartite graph with $k\geq3$. Then $smp(G)=2$ and each optimal strong matching preclusion set of $G$ is a set of two vertices from the same  bipartition set.}

\vskip.2cm
By Definition 2.1 and Lemma 2.6, we have the following theorem.
\vskip.2cm
{\bf Theorem 2.7} {\em $smp(BS_n)=2$ for $n\geq3$ and each optimal strong matching preclusion set of $BS_n$ is a set of two vertices from the same  bipartition set.}

\vskip.3cm

\n{\large\bf 3.\quad Matching preclusion of $BS_n$}

\vskip.2cm

In this section, we will consider the matching preclusion number and the optimal matching
preclusion set of $BS_n$.

Let $\Gamma$ be a finite group and $\Lambda\subseteq \Gamma$ such that the identify of $\Gamma$ is in $\Gamma-\Lambda$. The digraph {\em $Cay(\Gamma,\Lambda)$} is a Cayley graph with vertex set $\Gamma$ and arc set $\{\langle g,g\circ \lambda\rangle~|~g\in \Gamma,~\lambda\in \Lambda\}$. Let $\Lambda^{-1}=\{\lambda^{-1}~|~\lambda\in \Lambda\}$. The graph $Cay(\Gamma,\Lambda)$ is an undirected graph in the case that $\Lambda=\Lambda^{-1}$.

Let $S$ be the group of all permutations on $[1,n]$ with $n\geq2$ and $T$ be a set of transpositions of $S$. Let $G(T)$ be a graph with vertex set $\{u_1,u_2,\cdots,u_n\}$ and edge set $\{(u_i,u_j)~|~\langle i,j\rangle\in T\}$ [15]. The graph $G(T)$ is the {\em transposition generating graph} of the Cayley graph $Cay(S,T)$. By Definition 2.1 and Lemma 2.3, the $n$-dimensional bubble-sort star graph $BS_n$ is a Cayley graph with transposition generating graph $G(T)\cong K_{1,n-1}\cup P_n$ (see Fig. 2), where $T=\{\langle 1,i_1\rangle,\langle i_2,i_2+1\rangle~|~i_1\in[2,n],i_2\in[2,n-1]\}$.
\vskip.2cm

\begin{figure}[ht]
   \begin{center}
    \includegraphics[scale=1]{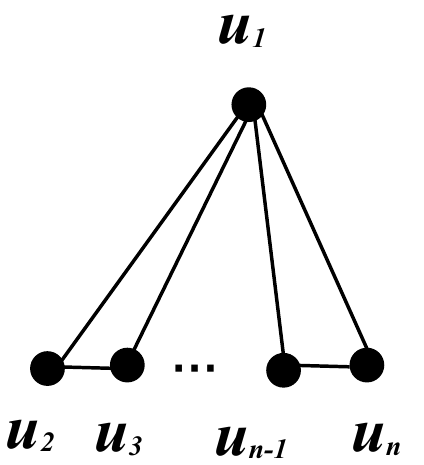}
    \end{center}
    \vspace{-0.5cm}\caption{\label{F2-5} $K_{1,n-1}\cup P_n$.}
\end{figure}

\vskip.2cm
{\bf Lemma 3.1 [16]} {\em  Let $G$ be a Cayley graph, which can be obtained from a connected transposition generating graph on $\left[1,n\right]$ with $n\geq4$. Then for any two distinct vertices $u$ and $v$ in distinct two bipartition sets of $G$, there exists a Hamiltonian path connecting them.}
\vskip.2cm

The graph $K_{1,n-1}\cup P_n$ is connected and $BS_n$ can be obtained from the transposition generating graph $K_{1,n-1}\cup P_n$ for $n\geq4$. Hence by Lemma 3.1, we get that for any two distinct vertices $u,v$ in distinct bipartition sets of $BS_n$ ($n\geq4$), there is a Hamiltonian path between $u$ and $v$.

\vskip.2cm
{\bf Lemma 3.2 [16]} {\em  Let $G$ be a bipartite and $k$-regular graph with $k\geq3$. Then $mp(G)=k$.}
\vskip.2cm

By Definition 2.1 and Lemma 3.2, we immediately have the following theorem.
\vskip.2cm
{\bf Theorem 3.3} {\em $mp(BS_n)=2n-3$ for $n\geq3$. }

\vskip.2cm

{\bf Lemma 3.4} {\em  Every optimal matching preclusion set of $BS_3$ is trivial. }
\vskip.2cm
\p Let $F\subseteq E(BS_3)$ be an arbitrary optimal matching preclusion set of $BS_3$. By Theorem 3.3, $|F|=mp(BS_3)=3$.
Suppose that $F$ is not trivial.
Let edges $a=(123,132)$, $b=(132,312)$, $c=(312,321)$, $d=(321,231)$, $e=(231,213)$, $f=(123,213)$, $g=(123,321)$, $h=(132,231)$ and $p=(213,312)$ (see Fig. 3). Now we will consider the following four cases.

\begin{figure}[ht]
   \begin{center}
    \includegraphics[scale=1]{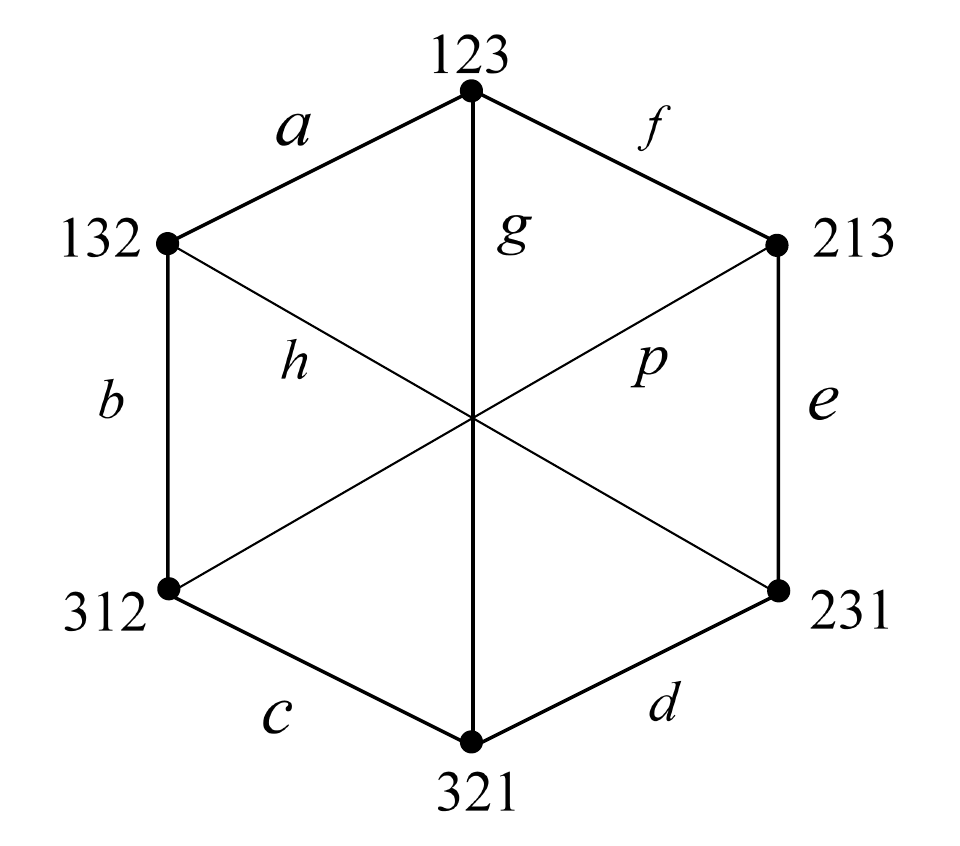}
    \end{center}
    \vspace{-0.5cm}\caption{\label{F2-5} Illustration of Lemma 3.4.}
\end{figure}

\vskip.2cm

{\bf Case 1.} {\em $|F\cap\{a,b,c,d,e,f\}|=3$.}

In this case, $\left\{g,h,p\right\}$ is a perfect matching in $BS_3-F$, a contradiction.

{\bf Case 2.} {\em $|F\cap\{g,h,p\}|=3$.}

In this case, both $\left\{a,c,e\right\}$ and $\left\{b,d,f\right\}$ are perfect matchings in $BS_3-F$, a contradiction.

{\bf Case 3.} {\em $|F\cap\{a,b,c,d,e,f\}|=2$ and $|F\cap\{g,h,p\}|=1$.}

In this case, the perfect matchings in $BS_3-F$ are listed in Table 1, a contradiction.
\vskip.2cm
\begin{center}
\setlength{\tabcolsep}{3.5pt}\renewcommand\arraystretch{1}
\begin{tabular}{cllllllllllllll}\hline

\multicolumn{1}{c}{$F$}&&\multicolumn{1}{c}{Perfect matching}&\multicolumn{11}{l}{\scriptsize }
\\\hline
\multicolumn{1}{c}{$\{a,b\}\subseteq F$ and $|F\cap \{g,p\}|=1$}&&\multicolumn{1}{c}{$\left\{f,h,c\right\}$}
\\
\multicolumn{1}{c}{$\{a,c\}\subseteq F$ and $|F\cap \{g,h,p\}|=1$}&&\multicolumn{1}{c}{$\left\{b,d,f\right\}$}
\\
\multicolumn{1}{c}{$\{a,d\}\subseteq F$ and $|F\cap\{g,h,p\}|=1$}&&\multicolumn{1}{c}{$\left\{f,h,c\right\}$ or $\left\{b,g,e\right\}$}
\\
\multicolumn{1}{c}{$\{a,e\}\subseteq F$ and $|F\cap\left\{g,h,p\right\}|=1$}&&\multicolumn{1}{c}{$\left\{b,d,f\right\}$}
\\
\multicolumn{1}{c}{$\left\{a,f\right\}\subseteq F$ and $|F\cap\left\{h,p\right\}|=1$}&&\multicolumn{1}{c}{$\left\{b,g,e\right\}$}
\\\hline
\end{tabular}

\vskip 0.2cm
Table 1. Perfect matchings in $BS_3-F$ in Case 3.
\end{center}

{\bf Case 4.}  {\em $|F\cap\{a,b,c,d,e,f\}|=1$ and $|F\cap\{g,h,p\}|=2$.}

In this case, $\left\{a,c,e\right\}$ or $\left\{b,d,f\right\}$ is a perfect matching in $BS_3-F$, a contradiction.

Hence every optimal matching preclusion set of $BS_3$ is trivial. \q

\vskip.2cm

{\bf Lemma 3.5} {\em Every optimal matching preclusion set of $BS_4$ is trivial.}

\vskip.2cm

\p Let $F\subseteq E(BS_4)$ be an arbitrary optimal matching preclusion set of $BS_4$. By Theorem 3.3, $\left|F\right|=mp(BS_4)=5$. 
Let $M^+=\{(v,v\circ\langle1,4\rangle)~|~v\in V(BS_4)\}$ and $M^-=\{(v,v\circ\langle3,4\rangle)~|~v\in V(BS_4)\}$. It is clearly that both $M^+$ and $M^-$ are perfect matchings of $BS_4$. Since there is no perfect matchings in $BS_4-F$, $|F\cap M^+|\geq1$ and $|F\cap M^-|\geq1$. By Lemma 2.5, $(F\cap M^+)\cap(F\cap M^-)=\emptyset$.
If $BS_4^i-F$ has a perfect matching $M_i$ for every $i\in[1,4]$, then $M=M_1\cup M_2\cup M_3\cup M_4$ is a perfect matching in $BS_4-F$, a contradiction. Hence there exists an integer $j\in[1,4]$, such that $BS_4^j-F$  has no perfect matchings. Without loss of generality, we assume $j=1$. By Theorem 3.3, $|F\cap E(BS_4^1)|\geq3$. Since $|F|=5$, $|F\cap M^+|\geq1$, $|F\cap M^-|\geq1$ and $(F\cap M^+)\cap(F\cap M^-)=\emptyset$, we have $|F\cap E(BS_4^1)|=3$, $|F\cap M^+|=|F\cap M^-|=1$. Hence $F\cap E(BS_4^1)$ is an optimal matching preclusion set of $BS_4^1$ and there is an isolated vertex $u$ in $BS_4^1-F$ by Lemma 3.4.

If $F\cap M^+=\{(u,u\circ\langle1,4\rangle)\}$ and $F\cap M^-=\{(u,u\circ\langle3,4\rangle)\}$, then $u$ is an isolated vertex in $BS_4-F$. Thus $BS_4-F$ does not have perfect matchings and $F$ is trivial. Now we suppose that $F$ is not trivial.

Since $BS_4$ is vertex transitive by Lemma 2.2, we assume $u=4321$. Let $u_1=4312$, $u_2=1324$, $v_1=3421$, $v_2=4231$, $v_3=2341$, $v_1^\prime=3412$, $v_2^\prime=1234$, $v_3^\prime=1342$ and $v_3^{\prime\prime}=2314$ (see Fig. 4). Then $F\cap E(BS_4^1)=\{(u,v_1),(u,v_2),(u,v_3)\}$, $\{(u,u_1),(v_1,v_1^{'}),(v_3,v_3^{''}),(v_3^\prime,u_2)\}\subseteq M^-$ and $\{(u,u_2),(v_2,v_2^{'}),(v_3,v_3^{'}),(u_1,v_3^{\prime\prime})\}\subseteq M^+$. Since $E(BS_4^{i})\cap F=\emptyset$ for every $i\in[2,4]$, there is a perfect matching $M_i$ in $BS_4^i-F$ by Theorem 3.3. Now we consider the following two cases.

\begin{figure}[ht]
   \begin{center}
    \includegraphics[scale=1]{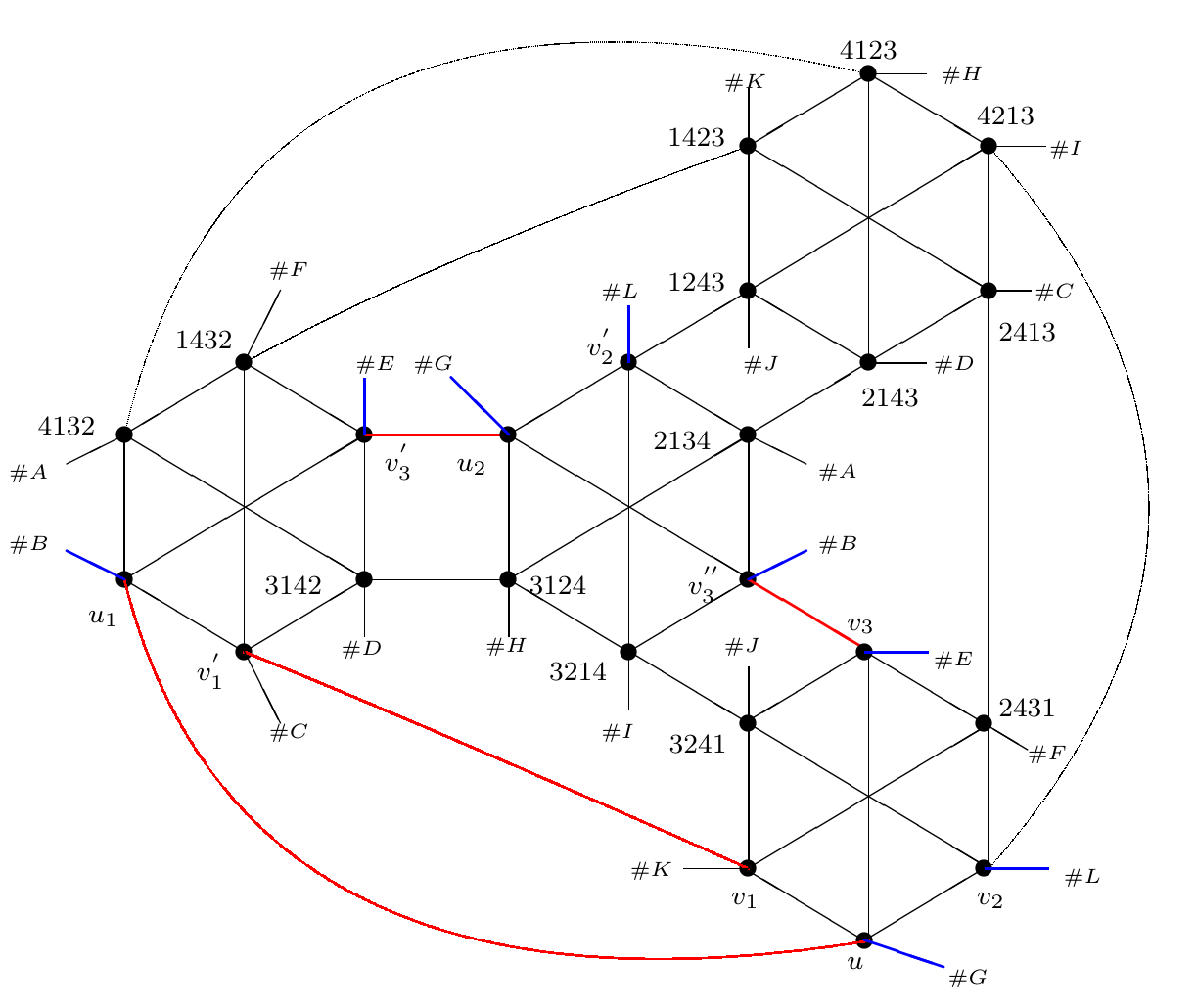}
    \end{center}
    \vspace{-0.5cm}\caption{\label{F2-5} Illustration of Lemma 3.5.}
\end{figure}

{\bf Case 1.} {\em $(u,u_1)\not\in F$.}

Suppose that $(v_1^{'},v_1)\not\in F$. Then $M_{12}=\{(u,u_1),(v_1^{'},v_1),(4132,1432),(3142,v_3^\prime),
(3241,v_3), \\(v_2,2431)\}$ is a perfect matching of $BS_4^{[1,2]}-F$. Thus $M_{12}\cup M_3\cup M_4$ is a perfect matching in $BS_4-F$, a contradiction. Hence $(v_1^{'},v_1)\in F$. Since $|F\cap M^-|=1$, $(v_3,v_3^{''})\not \in F$ and $(v_3^\prime,u_2)\not\in F$. Then $M_{124}=\{(u,u_1),(v_3,v_3^{''}),(v_3^\prime,u_2),(3241,v_1),(2431,v_2),(4132,1432),\\(v_1^{'},3142),(3124,3214),(v_2^{'},2134)\}$ is a perfect matching in $BS_4^{\{1,2,4\}}-F$. Thus $M_{124}\cup M_3$ is a perfect matching in $BS_4-F$,
a contradiction.

{\bf Case 2.} {\em $(u,u_1)\in F$.}

Since $F$ is not trivial, $(u,u_2)\not\in F$.
If $(v_2^{'},v_2)\not\in F$, then  $M_{14}=\{(u,u_2),(v_2^{'},v_2),\\(3124,3214),(2134,v_3^{\prime\prime}),(3241,v_1),(v_3,2431)\}$ is a perfect matching in $BS_4^{\{1,4\}}-F$. Thus $M_{14}\cup M_2\cup M_3$ is a perfect matching in $BS_4-F$, a contradiction. Hence $(v_2^{'},v_2)\in F$. Since $|F\cap M^+|=1$, $(v_3,v_3^{'})\not \in F$ and $(u_1,v_3^{\prime\prime})\not\in F$. Then $M_{124}'=\{(u,u_2),(v_3,v_3^{'}),(u_1,v_3^{\prime\prime}),(3241,v_1),(2431,v_2),(4132,1432),\\(v_1^{'},3142),(3124,3214),(v_2^{'},2134)\}$ is a perfect matching in $BS_4^{\{1,2,4\}}-F$. Thus $M_{124}'\cup M_3$ is a perfect matching in $BS_4-F$,
a contradiction.

Hence $(u,u_1)\in F$, $(u,u_2)\in F$ and $F$ is trivial.\q
\vskip.2cm

{\bf Theorem 3.6} {\em Every optimal matching preclusion set of $BS_n$ is trivial for $n\geq3$.}
\vskip.2cm

\p We prove this theorem by induction on $n$. For $n=3,4$, the theorem holds by Lemmas 3.4 and 3.5. Now we assume that $n\geq5$.
Let $F\subseteq E(BS_n)$ be an arbitrary optimal matching preclusion set of $BS_n$. By Theorem 3.3, $|F|=mp(BS_n)=2n-3$. 
Let $M^+=\{(v,v\circ\langle1,n\rangle)~|~v\in V(BS_n)\}$ and $M^-=\{(v,v\circ\langle n-1,n\rangle)~|~v\in V(BS_n)\}$. It is clearly that both $M^+$ and $M^-$ are perfect matchings of $BS_n$. Since there is no perfect matchings in $BS_n-F$, $|F\cap M^+|\geq1$ and $|F\cap M^-|\geq1$. By Lemma 2.5, $(F\cap M^+)\cap(F\cap M^-)=\emptyset$.

If $BS_n^{i}-F$ has a perfect matching $M_i$ for every $i\in[1,n]$, then $BS_n-F$ has a perfect matching $M_1\cup M_2\cup\cdots \cup M_n$, a contradiction. Thus there exists an integer $j\in[1,n]$ such that $BS_n^j-F$ does not have perfect matchings. Without loss of generality, we assume $j=1$. By Theorem 3.3, $|F\cap E(BS_n^1)|\geq2n-5$. Since $|F|=2n-3$, $|F\cap M^+|\geq1$, $|F\cap M^-|\geq1$ and $(F\cap M^+)\cap(F\cap M^-)=\emptyset$, we have $|F\cap E(BS_n^1)|=2n-5$, $|F\cap M^+|=|F\cap M^-|=1$ and $F\cap E(BS_n^i)=\emptyset$ for every $i\in[2,n]$. Hence $F\cap E(BS_n^1)$ is an optimal matching preclusion set of $BS_n^1$ and there is an isolated vertex $u$ in $BS_n^1-F$ by induction hypothesis. If $F\cap M^+=\{(u,u\circ\langle1,n\rangle)\}$ and $F\cap M^-=\{(u,u\circ\langle n-1,n\rangle)\}$, then $u$ is an isolated vertex in $BS_n-F$.
Thus $BS_n-F$ does not have perfect matchings and $F$ is trivial. Now we suppose that $F$ is not trivial.

Since $BS_n$ is vertex transitive by Lemma 2.2, we assume $u=n(n-1)\cdots321$.
Since $F$ is not trivial, there exists a vertex $v_2\in\{u\circ\langle1,n\rangle,u\circ\langle n-1,n\rangle\}$ such that $(u,v_2)\not\in F$. Without loss of generality, we assume $v_2\in V(BS_n^2)$. Since $BS_n$ is a bipartite graph, there exist two bipartition subsets $V_1,V_2\subseteq V(BS_n)$ with $u\in V_1$, $v_2\in V_2$ such that $V_1\cap V_2=\emptyset$, $V_1\cup V_2=V(BS_n)$ and every edge in $BS_n$ has one end vertex in $V_1$ and the other in $V_2$.

\vskip.2cm

\begin{figure}[ht]
   \begin{center}
    \includegraphics[scale=1]{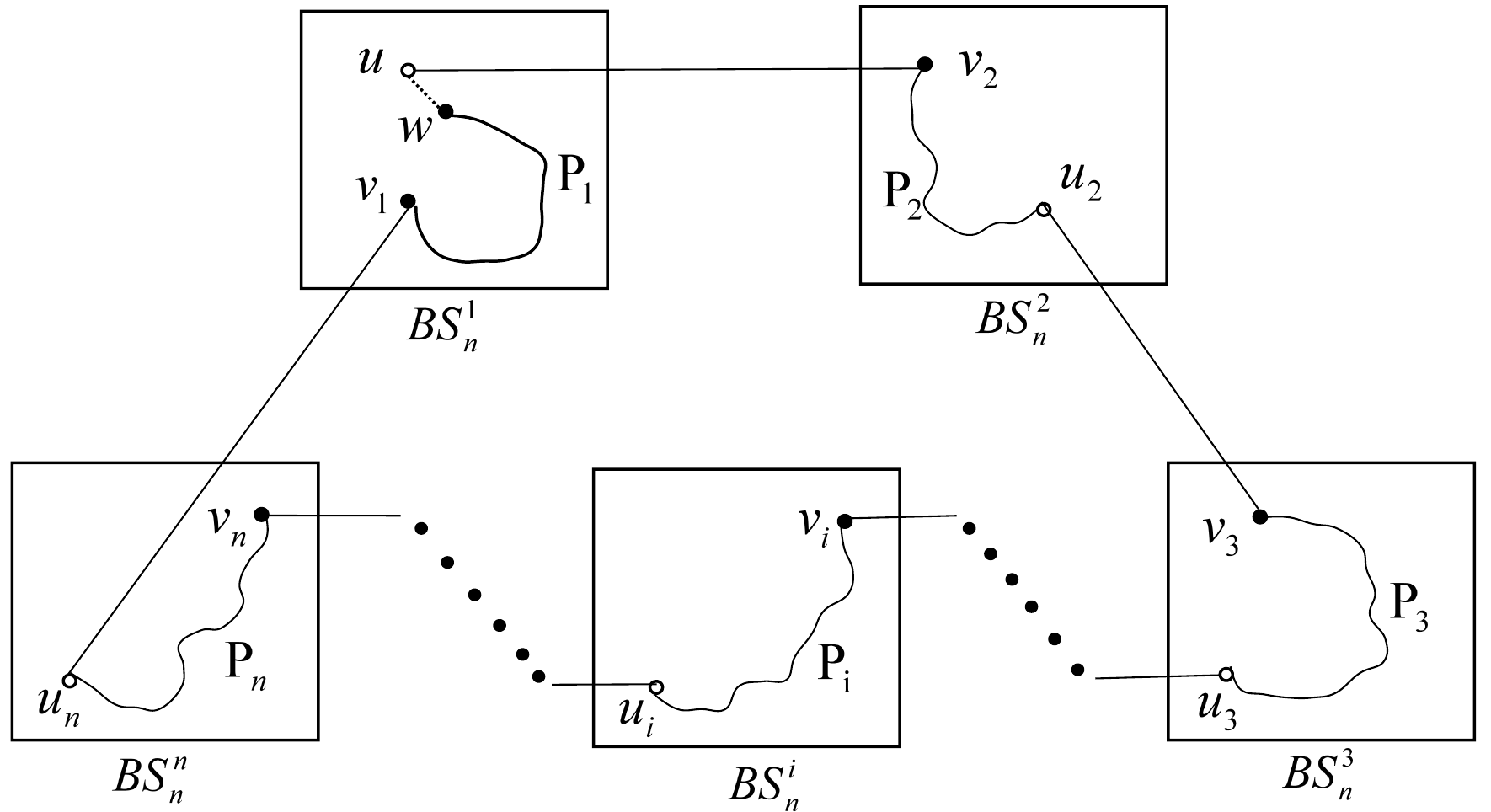}
    \end{center}
    \vspace{-0.5cm}\caption{\label{F2-5} A Hamiltonian path in $BS_n-F$.}
\end{figure}

\vskip.2cm

Let $E_{i,j}=\{(w_1,w_2)\in E(BS_n)~|~w_1\in V(BS_n^i),w_2\in V(BS_n^j)\}$ for $i,j\in[1,n]$ with $i\ne j$. By Lemma 2.4, $|E_{i,j}|=2(n-2)!$ for every $i,j\in[1,n]$ with $i\ne j$. By Definition 2.1, $|M^+\cap E_{i,j}|=|M^-\cap E_{i,j}|=(n-2)!$ for every $i,j\in[1,n]$ with $i\ne j$. Since $|F\cap M^+|=|F\cap M^-|=1$ and $(n-2)!>1$ for $n\geq5$, there exist vertices $u_k\in V(BS_n^k)\cap V_1$ and $v_{k+1}\in V(BS_n^{k+1})\cap V_2$ such that $(u_k,v_{k+1})\in M^+\cup M^--F$ for every $k\in[2,n-1]$. Also we can get two vertices $u_n\in V(BS_n^n)\cap V_1$ and $v_1\in V(BS_n^1)\cap V_2$ such that $(u_n,v_1)\in M^+\cup M^--F$ (see Fig. 5).

By Lemmas 2.3 and 3.1, there exists a Hamiltonian path $P_k$ between $v_k$ and $u_k$ in $BS_n^k$ for every $k\in[2,n]$. Also there exists a Hamiltonian path $P_1$ between $v_1$ and $u$ in $BS_n^1$. Let $w$ be the vertex such that $(u,w)\in E(P_1)$. Since $F\cap E(BS_n^1)$ is an optimal matching preclusion set of $BS_n^1$, $(E(P_1)-\{(u,w)\})\cap F=\emptyset$. Since $F\cap E(BS_n^k)=\emptyset$ for every $k\in[2,n]$, $E(P_k)\cap F=\emptyset$. Let $H$ be the subgraph of $BS_n$ with vertex set $V(H)=V(BS_n)$ and edge set $E(H)=\cup_{k=2}^{n-1} (E(P_i)\cup\{(u_i,v_{i+1})\})\cup E(P_1)\cup E(P_n)\cup\{(v_1,u_n),(u,v_2)\}-\{(u,w)\}$. Then $H$ is a Hamiltonian path between $u$ and $w$ in $BS_n-F$. Since $|V(H)|=|V(BS_n)|=n!$ for $n\geq5$, there is a perfect matching in $H$, which is also a perfect matching in $BS_n-F$, a contradiction.

Hence every optimal matching preclusion set of $BS_n$ is trivial.\q

\vskip.3cm

\n{\large\bf 4.\quad Conclusion}
\vskip.2cm
In this paper, we study the strong matching preclusion and the matching preclusion of the $n$-dimensional bubble-sort star graph $BS_n$. We show that the strong matching preclusion number of $BS_n$ is $2$ for $n\geq3$ and each optimal strong matching
preclusion set of $BS_n$ is a set of two vertices from the same bipartition set.
We also show that the matching preclusion number of $BS_n$ is $2n-3$ for $n\geq3$ and that every optimal matching preclusion set of $BS_n$ is trivial. The (strong) matching preclusion of interconnection network determines the fault tolerance and stability of the network.

Our findings are meaningful for establishing and optimizing the network topology of the distributed computer system, and also help to improve the fault tolerance and stability of the network. In the follow-up, we plan to further study the theory and application of the matching preclusion in graph theory to optimize the distributed networks.

\n{\large\bf Acknowledgements}
\vskip.2cm
This work is supported by NNSF of China (Nos. 11801450), Natural Science Foundation of Shaanxi Province, China (Nos. 2019JQ-506).

\end{document}